\documentclass[a4paper,12pt,reqno]{article}
\usepackage{amsthm}
\usepackage{amsmath,amssymb,latexsym,amsfonts,mathrsfs}
\usepackage{bm}
\usepackage{bigints}
\usepackage{color}
\usepackage{ulem}
\usepackage{mathtools}
\usepackage[top=30truemm,bottom=30truemm,left=20truemm,right=20truemm]{geometry}
\numberwithin{equation}{section}

\newcommand{\R}{\mathbb{R}}

\newcommand{\ep}{\varepsilon}

\theoremstyle{plain}
\newtheorem{Thm}{{\bf Theorem}}[section]

\newtheorem{Lem}[Thm]{{{\bf Lemma}}}
\newtheorem{Prop}[Thm]{{\bf Proposition}}

\newtheorem{Ass}[Thm]{{\bf Assumption}}

\newtheorem{Rem}[Thm]{{\bf Remark}}

\newcounter{Exami}

\renewcommand{\Im}{\mathrm{Im}}
\renewcommand{\Re}{\mathrm{Re}}
\allowdisplaybreaks[2]
\begin{document}
\fontencoding{T1}\selectfont
\title{Non-smoothness of the fundamental solutions for Schr\"{o}dinger equations with
super-quadratic and spherically symmetric potential
}
\author{Keiichi Kato\thanks{Department of Mathematics, Faculty of Science, Tokyo University of Science, Kagurazaka 1-3,
Shinjuku-ku, Tokyo 162-8601, Japan. E-mail: kato@rs.tus.ac.jp}, 
Wataru Nakahashi\thanks{Department of Mathematics, Graduate School of Science, Tokyo University of Science, Kagurazaka 1-3,
Shinjuku-ku, Tokyo 162-8601, Japan. Email: 1122705@ed.tus.ac.jp} 
 and Yukihide Tadano\thanks{
Graduate School of Science, University of Hyogo, Shosha,
Himeji, Hyogo 671-2280, Japan. E-mail: tadano@sci.u-hyogo.ac.jp}}
\maketitle
%
\begin{abstract} We study non-smoothness of the fundamental solution  for the Schr\"{o}dinger equation with a spherically symmetric and super-quadratic potential in the sense that $V(x)\geq C|x|^{2+\varepsilon}$ at infinity with constants $C>0 $ and $\varepsilon>0$. More precisely, we show the fundamental solution $E(t,x,y)$ does not belong to $C^{1}$ as a function of $(t,x,y)$, which partially solves Yajima's conjecture.
\end{abstract}
\section{Introduction}
\par We consider the following initial value problem for the Schr\"{o}dinger equation on $\R^d$ with a super-quadratic and spherically symmetric potential $V(x)$
\begin{align}\label{eq1}
\begin{dcases}
i \frac{\partial u}{\partial t}(t,x)  = - \Delta u(t,x)+ V(x) u(t,x),\quad(t,x)\in \R \times \R^d,
\\u(0,\cdot)=u_{0} \in  L^2(\R^d),
\end{dcases}
\end{align}
where $d\geq3$, $i=\sqrt{-1}$, $\Delta =\sum\limits_{j=1}^{d}\partial^{2}/\partial{x_j}^{2}$, and $u:\mathbb R\times \mathbb R^{d} \rightarrow \mathbb C$.
The aim of this paper is to show non-smoothness of the fundamental solutions under the following assumption, which partially solves Yajima's conjecture \cite{Y1}.

\begin{Ass}\label{Ass}
The potential $V(x)\in C^3(\R^d)$ is a real valued spherically symmetric function.
Moreover, letting $V(x)=\widetilde{V}(|x|)$, there are constants $R>0$ and $c>1$ such that for $r \in [R,\infty)$
\begin{align}
&\widetilde{V}^{\prime\prime}(r)>0,\label{Ass1}
\\ 
&r\widetilde{V}^{\prime}(r)\geq 2c\widetilde{V}(r)>0,\label{Ass2}
\\
&\widetilde{V}^{(j)}(r)=O\left(\frac{1}{r}\right)\widetilde{V}^{(j-1)}(r) ,\ j=1,2,3 \label{Ass3}.
\end{align}
\end{Ass}
It follows from \eqref{Ass2} that there exists $C>0$ such that $V(x)\geq C|x|^{2c}$ for $|x|\gg1$, which implies $H=-\Delta+V$ is essentially self-adjoint on $C_0^\infty(\R^d)$.
Thus the equation generates a unique unitary propagator $U(t)=e^{-itH}$ on $L^{2}(\mathbb R^{d})$ and $u(t,x)=(U(t)u_{0})(x)$ is the unique solution of (\ref{eq1}).
We denote the integral kernel of $U(t)$ by $E(t,x,y)$, that is,
\[
(U(t)u_0)(x)= \int_{\R^d} E(t,x,y) u_0(y) dy,
\]
and we call $E(t,x,y)$ the fundamental solution of (\ref{eq1}). 
\par
The main theorem of this paper is the following, which claims that $E(t,x,y)$ is generically nowhere $C^1$.
\begin{Thm}\label{main}
Suppose that $d\geq3$ and that $V(x)$ satisfies Assumption \ref{Ass}.
Then, for any $t_{0} \in\mathbb R$ and $r_{1},r_{2}>0$, there exist $x_{0}, y_{0}\in \mathbb R^{d}$, satisfying $|x_{0}|=r_{1}$ and $|y_{0}|=r_{2}$, such that the fundamental solution of (\ref{eq1}) does not belong to $C^1$ near $(t_{0}, x_{0}, y_{0})$. 
\end{Thm}
Moreover, we can obtain the following proposition which is (almost) stronger than Theorem~\ref{main}.
\begin{Prop}\label{cor-main}
Suppose that $d\geq3$ and that $V(x)$ satisfies Assumption \ref{Ass}.
Let $Y_{n}$ be a spherical harmonic of arbitrary degree $n$ with $\|Y_{n}\|_{L^2(\mathbb{S}^{d-1})}=1$ and let
\begin{align*}
E^{Y_{n}}(t,x,y) =& (E(t,|x|\cdot,y),Y_{n})_{L^2(\mathbb{S}^{d-1})} Y_{n}\left(\frac{x}{|x|}\right) .
\end{align*}
Then, for any $(t_{0}, x_{0}, y_{0})\in\mathbb R\times\mathbb R^{d}\times\mathbb R^{d}$, $E^{Y_{n}}$ is not in $C^1$ near $(t_{0}, x_{0}, y_{0})$.
\end{Prop}

Here $(\cdot,\cdot)_\mathcal{H}$ denotes the inner product of a Hilbert space $\mathcal{H}$, and we suppose $(\cdot,\cdot)_\mathcal{H}$ is linear with respect to the first entry and anti-linear with respect to the second entry.

We give some comments on what Theorem~\ref{main} and Proposition~\ref{cor-main} mean and how our proof works.
We denote the self-adjoint extension of $H=-\Delta+V$ by the same symbol.
Since $V(x)\to\infty$ as $|x|\to\infty$ under Assumption \ref{Ass}, the spectrum $\sigma(H)$ of $H$ is discrete.
Thus we have, at least formally,
\begin{align*}
E(t,x,y)=\sum\limits_{l}e^{-i\lambda_{l} t}u_{l}(x)\overline{u_{l}(y)},
\end{align*}
where $\{u_l\}_l$ is an arbitrary complete orthonormal system of $L^2(\R^d)$ with each $u_l$ being an eigenfunction of $H$: $Hu_l=\lambda_l u_l$.
\par
It is known that $H$ is decomposed by the partial wave expansion; let $\{Y_{nm} \mid n=0,1,2,\dots,\ m=1,2,\dots,d_n\}$ be a complete orthonormal system of $L^2(\mathbb{S}^{d-1})$ such that each $Y_{nm}$ is a spherical harmonic of degree $n$, where $d_n=\binom{d+n-1}{d-1}-\binom{d+n-3}{d-1}$ is the dimension of the space of the spherical harmonics of degree $n$.
Then we have
\begin{align}
E(t,x,y)=\sum_{n=0}^\infty \sum_{m=1}^{d_n} \sum_{l=0}^\infty e^{-i\lambda_{n,l}t} |x|^{-\frac{d-1}{2}} f_{n,l}(|x|)Y_{nm}\left(\frac{x}{|x|}\right)\overline{|y|^{-\frac{d-1}{2}} f_{n,l}(|y|)Y_{nm}\left(\frac{y}{|y|}\right)} , \label{PWE}
\end{align}
where $f_{n,l}\in L^2((0,\infty))$ is the normalized $\lambda_{n,l}$-eigenfunction of the Schr\"odinger operator on the half line $(0,\infty)$:
\begin{equation}\label{HLS}
- f_{n,l}''(r) + \left(\frac{(d-1)(d-3)}{4r^{2}}+\frac{n(n+d-2)}{r^{2}}+\widetilde{V}(r)\right) f_{n,l}(r) = \lambda_{n,l} f_{n,l}(r), \quad r\in(0,\infty),
\end{equation}
with the boundary condition $f_{n,l}(0)=0$.
Hence, if $Y_n$ is a spherical harmonic of degree $n$, we have
\begin{align*}
E^{Y_n}(t,x,y) =& \sum_{l=0}^\infty e^{-i\lambda_{n,l}t} |x|^{-\frac{d-1}{2}} f_{n,l}(|x|) Y_n\left(\frac{x}{|x|}\right)\overline{|y|^{-\frac{d-1}{2}} f_{n,l}(|y|) Y_n\left(\frac{y}{|y|}\right)}\\
=& |x|^{-\frac{d-1}{2}} |y|^{-\frac{d-1}{2}} Y_n\left(\frac{x}{|x|}\right) \overline{Y_n\left(\frac{y}{|y|}\right)} \sum_{l=0}^\infty e^{-i\lambda_{n,l}t} f_{n,l}(|x|) \overline{f_{n,l}(|y|)}.
\end{align*}
Since $\sum_{l=0}^\infty e^{-i\lambda_{n,l}t} f_{n,l}(r) \overline{f_{n,l}(s)}$ is the fundamental solution associated to the above one-dimensional Schr\"odinger operator,
we can employ the argument of Yajima \cite{Y1}.

One might expect that Proposition~\ref{cor-main} and the representation \eqref{PWE} "imply" that our $E(t,x,y)$, the formal summation over spherical harmonics
\begin{align}
E(t,x,y)=\sum_{n=0}^\infty \sum_{m=1}^{d_n} E^{Y_{nm}}(t,x,y), \label{PWE2}
\end{align}
would be nowhere in $C^1$, which is, however, still not clear.
Instead, we can obtain Theorem~\ref{main}.

\begin{Rem}\label{rem-thm}
(i) $E^{Y_{n}}(t,x,y)$ is actually defined as a distribution on $\R^{2d+1}$ by
\begin{align*}
\langle E^{Y_{n}}, F \rangle 
=& \int_{\R} \int_{\R^d} \int_{\R^d} \int_{\mathbb{S}^{d-1}} E(t,|x|\omega,y) \overline{Y_{n}(\omega)} Y_{n}\left(\frac{x}{|x|}\right) F(t,x,y) d\omega dy dx dt \\
=& \int_{\R} \int_{\R^d} \int_{\mathbb{S}^{d-1}} \overline{Y_{n}(\omega)} Y_{n}\left(\frac{x}{|x|}\right) \left[U(t)\left(F(t,x,\cdot)\right)\right](|x|\omega) d\omega dx dt .
\end{align*}
\par
(ii) The representation \eqref{PWE} implies the identity
\begin{align*}
E^{Y_{n}}(t,x,y) =& (E(t,x,|y|\cdot),\overline{Y_{n}})_{L^2(\mathbb{S}^{d-1})} \overline{Y_{n}\left(\frac{y}{|y|}\right)} .
\end{align*}
\end{Rem}
%
%
%
%
%
%
%
%
%
\par
The smoothness of the fundamental solution is related to the growth rate of $V$.
Fujiwara \cite[Theorem 1.1]{F} has given the construction of the fundamental solution $E(t,x,y)$ with the classical orbit if $V$ is at most of quadratic growth, which shows that $E(t,x,y)$ is smooth with respect to $(x,y)$ for any $t\neq0$ small enough as a corollary. 
\par
We can give an alternative proof for the smoothness of the fundamental solution 
by \cite{IK}, which shows characterization of wave front set of solutions to Schr\"odinger equations with sub-quadratic potentials via wave packet transform. 
Theorem 1.2 in \cite{IK} yields that $( x_0, \xi_0)$ is not in the wave front set of $E(t_0,\cdot )$
if and only if 
$$
|e^{-it_0{\triangle}}\varphi_\lambda(x(0;t_0,x_0,\lambda \xi_0 ))|
\le C \lambda^{-N} \text{ for } \lambda \ge 1, 
$$
where $\varphi$ is a Schwartz function, $\varphi_\lambda (x)= \lambda^{d/2}\varphi(\lambda^{1/2}x)$ and $x(t;t_0, x_0, \lambda \xi_0 )$ is 
the solution of $\dot{x}(t)=\xi (t), \dot{\xi}(t)=-\nabla V(t,x) , x(t_0)=x_0, \xi(t_0)= \lambda \xi_0 $, 
which implies the smoothness of the fundamental solutions for Schr\"odinger equations with sub-quadratic potentials.  
\par
On the other hand, Yajima \cite[Theorem 1.2]{Y1} has studied that if $V$ is super-quadratic and the spatial dimension is one, then $E(t,x,y)$ is not smooth anywhere with respect to $(t,x,y)$. 
The proof is given by the estimates of eigenvalues and eigenfunctions of $H$. 
Yajima has conjectured that the same result as in Yajima \cite[Theorem 1.2]{Y1} is true even for higher dimensional cases.
\par
As a first step for generalization of Yajima's result, we treat the case that $V$ is spherically symmetric with the dimension $d\geq3$. 
We use the unitary equivalence of $H=-\Delta+V$ to $\bigoplus_{n,m} H_{nm}$, where the operator $H_{nm}$ is defined on $(0,\infty)$ (see Lemma \ref{lem-decomp}), and the estimates of eigenvalues and eigenfunctions which are shown by the same way as in \cite{Y1}. 
Our main theorem partially solves Yajima's conjecture.
In the case that $d=2$, the projection of $H$ onto the subspace of spherically symmetric functions is unitarily equivalent to $H_0=-\frac{d^2}{dr^2}-\frac{1}{4r^2}+\widetilde{V}(r)$ on $L^2((0,\infty))$ (see Section~{\ref{sec-polar}}).
The potential $-\frac{1}{4r^2}+\widetilde{V}(r)$ is not bounded from below near $0$, which requires additional argument of self-adjoint extensions.
We shall discuss the two-dimensional case in the forthcoming paper.
%
%
%
%
%
\par
For a Schr\"{o}dinger operator $-\Delta_{g}$ on a complete Riemannian manifold $(M,g)$, the smoothness of the fundamental solution depends on whether $M$ is compact or not.
Kapitanski and Rodnianski \cite[Theorem I-I\hspace{-1.2pt}I\hspace{-1.2pt}I]{KR}
has studied that $E(t,x,y)$ is not smooth if $M$ is the circle.
Taylor \cite[Section 1]{MT} has mentioned that $E(t,x,y)$ is not smooth if $M$ is the sphere.
Yajima \cite[Remark 4]{Y1} has pointed out that $E(t,x,y)$ is not smooth if $M$ is the bounded interval $[0,\pi]$ with the Dirichlet condition. 
Taira \cite[Remark 3.3]{KT} has studied that $E(t,x,y)$ is not smooth if $M$ is compact. 
When $M$ is non-compact,
Doi \cite[Theorem 1.5]{D} has studied the smoothness of the fundamental solution in terms of the wave front set.
Taira \cite[Theorems~1.1]{KT} has given a sufficient condition under which $E(t,x,y)$ is smooth.
\par
We introduce some notation. 
For sets $U$ and $V$, we write $U\Subset V$ if $U$ is relatively compact with respect to $V$.
We write $C_{0}^{\infty}(U)=\{f\in C^{\infty}(U)\mid \operatorname{supp}f \Subset  U\}$, where $\operatorname{supp} f$ denotes the support of $f$.
For an open interval $I$, we denote the Sobolev space on $I$ of order $k\in\mathbb{N}$ by $H^{k}(I)=\{f\in L^2(I)\mid f^{(m)} \in L^2(I) \text{ for } m=1,\dots,k\}$.
We denote by $H_{0}^{1}(I)$ the closure of $C_{0}^{\infty}(I)$ in $H^{1}(I)$. 
We denote the Schwartz space by $\mathcal S(\mathbb R^{d})$.
For any $f\in\mathcal{S}(\mathbb{R})$, the Fourier transform $\widehat{f}$ of $f$ is defined by $\widehat{f}(\lambda)=\int_{\R} e^{-it\lambda}f(t) dt$. 
\par
This paper is organized as follows.
In Section~\ref{sec-polar}, we introduce the polar coordinates to decompose the operator $H=-\Delta+V$ into the direct sum of operators on the half line.
In Section~\ref{sec-ev ef}, we observe estimates of eigenvalues and eigenfunctions. 
In Section~\ref{sec-proof}, we prove Theorem~\ref{main} and Proposition~\ref{cor-main}.
%
%
\section{Schr\"odinger operator in the polar coordinates}\label{sec-polar}
In the following, space dimension $d$ is larger than or equal to 3.
We call a function $Y$ on $\mathbb{S}^{d-1}$ \textit{spherical harmonic of degree n} if $Y$ is the restriction to $\mathbb{S}^{d-1}$ of a homogeneous harmonic polynomial of degree $n$. 
We denote by $\mathcal{H}_{n}$ the space of spherical harmonics of degree $n$. 
The following lemma is well-known (for more details and proofs, see Stein and Weiss \cite[Section 2 in Chapter I\hspace{-1.2pt}V]{SW} or Yajima \cite[Section 13 in Chapter 4]{Y2}). 
\renewcommand{\labelenumi}{(\theenumi)}
\begin{Lem}\label{lem HS}
Let $\Delta_{\mathbb S^{d-1}}$ be the Laplace-Beltrami operator on $\mathbb S^{d-1}$.
Then:
\renewcommand{\theenumi}{\roman{enumi}}
\begin{enumerate}
\item $-\Delta_{\mathbb S^{d-1}}Y=n(n+d-2)Y$ for $Y\in \mathcal{H}_{n}$.
\item $L^2(\mathbb{S}^{d-1}) = \bigoplus_{n=0}^\infty \mathcal{H}_{n}$.
\item $d_n = \dim \mathcal{H}_{n}=\binom{d+n-1}{d-1}-\binom{d+n-3}{d-1}$. 
\end{enumerate}
\end{Lem}
Let $\{Y_{nm} \mid m=1,\ldots, d_{n}\}$ be an orthogonal basis of $\mathcal{H}_{n}$, then $\{Y_{nm} \mid n=0,1,\ldots \text{ and } m=1,\ldots, d_{n}\}$ is a complete orthonormal system of $L^{2}(\mathbb S^{d-1})$.
If we define $J_{nm}:L^2(\R^d)\to L^2((0,\infty))$ by
\[
J_{nm}u(r)=r^{\frac{d-1}{2}}\left(u(r\cdot),Y_{nm}\right)_{L^{2}(\mathbb{S}^{d-1})},
\]
we easily see that the operator
\[
J:L^{2}(\mathbb{R}^{d})\ni u \longmapsto \left\{J_{nm}u(r)\right\}_{nm}\in\bigoplus_{n,m}L^2((0,\infty))
\]
is a unitary operator, and that the adjoint is given by
\[
J^*(\{g_{nm}\})(x)
=\sum_{n,m} J_{nm}^*(g_{nm})(x)
=\sum_{n,m} |x|^{-\frac{d-1}{2}} g_{nm}(|x|)Y_{nm}\left(\frac{x}{|x|}\right).
\]
Since $\Delta=\partial_{r}^{2}+\frac{d-1}{r}\partial_{r}+\frac{1}{r^{2}}\Delta_{\mathbb S^{d-1}}$ holds in the polar coordinates $(r,\omega)\in (0,\infty)\times\mathbb{S}^{d-1}$, we have
\[
J_{nm}H J_{nm}^* g(r) = - g''(r)+\left\{\frac{(d-1)(d-3)}{4r^{2}}+\frac{n(n+d-2)}{r^{2}}+\widetilde{V}(r)\right\}g(r).
\]
Since the above operators are independent of $m$, we let
\[
H_n = J_{nm}H J_{nm}^* = - \frac{d^2}{dr^2}+\frac{(d-1)(d-3)}{4r^{2}}+\frac{n(n+d-2)}{r^{2}}+\widetilde{V}(r).
\]
%
\begin{Lem}
If $(d,n)\neq(3,0)$, $H_n$ is essentially self-adjoint on $C_0^\infty((0,\infty))$. 
In particular, the domain of its closure is the maximal domain $\{g\in L^2((0,\infty)) \mid H_ng\in L^2((0,\infty))\}$.
If $(d,n)=(3,0)$, $H_0$ is essentially self-adjoint on $J_{00}C_0^\infty(\R^d)$, and the domain of its closure is $\{g\in H^{1}_{0}((0,\infty)) \mid H_ng\in L^2((0,\infty))\}$.
\end{Lem}

\begin{proof}
Weyl's limit point-limit circle criterion \cite[Theorem~X.7]{RS2} implies that $H_n$ is essentially self-adjoint on $C_0^\infty((0,\infty))$ if and only if $H_{n}$ is in the limit point case at both $0$ and $\infty$.
It follows from \cite[Theorem~X.8]{RS2} 
that $H_{n}$ is in the limit point case at $\infty$ for any $(d,n)$, since the potential of $H_{n}$ is bounded from below.
\cite[Theorem~X.10]{RS2} implies that $H_{n}$ is in the limit point case at $0$ if and only if the coefficient $\frac{(d-1)(d-3)}{4}+n(n+d-2)$ of $r^{-2}$ is not smaller than $\frac34$, i.e., $(d,n)\neq(3,0)$.
Thus we have the assertion in the case $(d,n)\neq(3,0)$.

On the other hand, if $(d,n)=(3,0)$, $H_{n}$ is in the limit point case at $\infty$ and is in the limit circle case at $0$.
In this case the self-adjoint extensions of $H_0=-\frac{d^2}{dr^2}+\widetilde{V}(r)$ on $C_0^\infty((0,\infty))$ are those with the domain
\[
\mathcal{D}_\theta=\{g\in H^1((0,\infty))\mid -g''+\widetilde{V}g\in L^2((0,\infty)), g(0)\cos\theta+g'(0)\sin\theta=0\},
\]
where $\theta\in[0,2\pi)$.
We can see that $C_0^\infty((0,\infty))\subset J_{00}C_0^\infty(\R^d)$ and that $J_{00}C_0^\infty(\R^d) \subset \mathcal{D}_\theta$ only if $\theta=0$, since $J_{00}u(0)=0$ and $(J_{00}u)'(0)=1$ if $u\in C_0^\infty(\R^d)$ equals to $|\mathbb{S}^{2}|^{-1/2}=(4\pi)^{-1/2}$ near $0$.
Therefore, the self-adjoint extension of $H_0$ on $J_{00}C_0^\infty(\R^d)$ is unique and the domain of its closure is $\mathcal{D}_0$.
\end{proof}

We denote the above-mentioned self-adjoint extensions by the same symbols.
\begin{Lem}\label{lem-decomp}
One has
\begin{align*}
JHJ^{*} = \bigoplus_{n,m} H_{nm},\quad \text{where } H_{nm}=H_{n}.
\end{align*}
More precisely, $J^*(\{g_{nm}\}_{nm})\in D(H)$ if and only if $g_{nm}\in D(H_{n})$ for any $n,m$ and $\sum\limits_{n,m}\|H_{n}g_{nm}\|^{2}<\infty$, and
 \begin{align*}
JHJ^*(\{g_{nm}\}_{nm}) = \{H_{n}g_{nm}\}_{nm},\quad J^*(\{g_{nm}\}_{nm})\in D(H).
\end{align*}
%
\end{Lem}
%
%
\section{Estimates of eigenvalues and eigenfunctions}\label{sec-ev ef}
In this section, we fix $n$, and consider the eigenvalue problem
\begin{align}\label{eigeneq}
- \frac{d^{2}}{dr^{2}}f(r)+\left\{\frac{(d-1)(d-3)}{4r^{2}}+\frac{n(n+d-2)}{r^{2}}+\widetilde{V}(r)\right\}f(r)=\lambda f(r)
\end{align}
with the boundary condition associated to the domain of $H_n$. The following argument is based on \cite[$\S 3$]{Y1}.
%
\subsection{Asymptotic behavior in bounded regions}
We first study the asymptotic behavior of the solutions of \eqref{eigeneq} in the subset 
\[
\Omega_{\lambda}=\left\{r\in[\lambda^{-1/4},\infty)\left|\frac{(d-1)(d-3)}{4r^{2}}+\frac{n(n+d-2)}{r^{2}}+\widetilde{V}(r)\leq\frac{\lambda}{2}\right.\right\}
\]
as $\lambda\to\infty$.
Let $\lambda_0$ be sufficiently large so that $\Omega_\lambda$ is an interval for $\lambda\geq\lambda_0$.
We write 
\[
U_{d,n}(r)=\frac{(d-1)(d-3)}{4r^{2}}+\frac{n(n+d-2)}{r^{2}}+\widetilde{V}(r)
\]
for simplicity. 
We set 
\begin{align*}
S(r)=\int_{1}^{r}\sqrt{\lambda-U_{d,n}(s)}ds, \quad
a(r)=\left(\lambda-U_{d,n}(r)\right)^{-1/4},
\end{align*}
and make change of variables $y=S(r)$ and $f(r)=a(r)w(y)$.
Then \eqref{eigeneq} implies
\begin{align}\label{eq4}
\frac{d^{2}w}{dy^{2}}+w+
\mathcal{U}
(y,\lambda)w=0,
\end{align}
where
\[
\mathcal{U}
(y,\lambda)=\frac{U^{\prime\prime}_{d,n}(r(y))}{4(\lambda-U_{d,n}(r(y)))^{2}}+\frac{5(U^{\prime}_{d,n}(r(y)))^{2}}{16(\lambda-U_{d,n}(r(y)))^{3}} .
\]
It follows that $w$ satisfies the integral equation
\begin{align}\label{inteq}
w(y)=w(0)\cos y+w^{\prime}(0)\sin y-\int_{0}^{y}\sin(y-z)\mathcal{U}(z,\lambda)w(z)dz.
\end{align}
Note that $w(0)$ and $w'(0)$ are given by $f(1)$, $f'(1)$, $U_{d,n}(1)$, $U_{d,n}'(1)$ and $\lambda$; in fact,
\begin{align*}
w(0) =& \left(\lambda-U_{d,n}(1)\right)^{1/4}f(1), \\
w'(0) =& \left.\frac{1}{S'(r)}\frac{d}{dr}\left(\left(\lambda-U_{d,n}(r)\right)^{1/4}f(r)\right)\right|_{r=1}.
\end{align*} 
We set $C_{\lambda}=(w(0)-iw^{\prime}(0))$ so that $w(0)\cos y+w^{\prime}(0)\sin y=\text{Re}(C_{\lambda}e^{iy})$.
%
\begin{Lem}\label{lem-cls-allowed}
Suppose that a real valued function $f$ satisfies \eqref{eigeneq} on $\Omega_\lambda$.
Then,
for any compact set $I$ in $\Omega_{\lambda}$ which is independent of $\lambda$,
\begin{align*}
f(r)=\lambda^{-1/4}\left(\Re\{C_{\lambda}e^{i\sqrt{\lambda}r}\}+O(C_{\lambda}\lambda^{-1/2})\right), \quad r\in I, \ \lambda\geq\lambda_0.
\end{align*}
\end{Lem}
\begin{proof}
The assertion is proved similarly to \cite[Lemma 3.1]{Y1}.
Let $R$ be the same as in Assumption~\ref{Ass}.
Since $U_{d,n}^{\prime}(r) \leq CU_{d,n}(r) \leq C(\lambda-U_{d,n}(r))$ on $\Omega_{\lambda}\cap[R,\infty)$, 
we have, for any $r_{0}\in\Omega_{\lambda}\cap(R,\infty)$,
\[
\int_{S(R)}^{S(r_{0})} |\mathcal{U}(z,\lambda)| dz \leq \frac{C}{\lambda^{3/2}}\int_{R}^{r_{0}}U_{d,n}^{\prime\prime}(r)dr+\int_{R}^{r_{0}}\frac{CU_{d,n}^{\prime}(r)}{\left(\lambda-U_{d,n}(r)\right)^{3/2}}dr \leq \frac{C}{\lambda^{1/2}}.
\]
On the other hand, for any $r_{1}\in\Omega_{\lambda}\cap(0,1)$, we have
\[
\int^{S(1)}_{S(r_{1})} |\mathcal{U}(z,\lambda)| dz \leq \frac{C}{\lambda^{3/2}}\int_{\lambda^{-1/4}}^{1}|U_{d,n}^{\prime\prime}(r)|dr+\frac{1}{\lambda^{5/2}}\int_{\lambda^{-1/4}}^{1}C|U_{d,n}^{\prime}(r)|^{2}dr \leq \frac{C}{\lambda^{1/2}}.
\]
Applying to (\ref{inteq}), there exists a constant $C>0$ independent of $\lambda\geq\lambda_{0}$ and $r\in\Omega_{\lambda}$ such that 
\begin{align*}
\left|a(r)^{-1}f(r)-\text{Re}(C_{\lambda}e^{iS(r)})\right|\leq C\left|C_{\lambda}\right|\lambda^{-1/2}.
\end{align*}
We have $S(r)=\sqrt{\lambda}r+O(\lambda^{-1/2})$ and $a(r)=\lambda^{-1/4}(1+O(\lambda^{-1}))$ uniformly on any compact set $I$, which implies the assertion.
\end{proof}
The next lemma gives a lower bound of $|C_{\lambda}|$ in the previous lemma. The assertion is proved by the same way as in \cite[Lemma 3.2]{Y1} and we omit the proof.
\begin{Lem}\label{lem-cls-allowed-const}
Suppose that a real valued function $f\in L^{2}((0,\infty))$ satisfies \eqref{eigeneq} and $\|f\|_{L^{2}((0,\infty))}=1$,
and let $C_\lambda$ be as in Lemma \ref{lem-cls-allowed}.
Then there is a constant $C>0$ such that
\begin{align*}
|C_{\lambda}|\geq C\lambda^{(c-1)/4c} , \quad \lambda\geq\lambda_{0},
\end{align*} 
where $c>1$ is the constant in Assumption \ref{Ass}.
\end{Lem}
%
%
%
\subsection{Gaps of the eigenvalues}
Let $\lambda_{n,0}<\lambda_{n,1}<\cdots$ be the eigenvalues of $H_n$.
We show that the gap $\lambda_{n,l+1}-\lambda_{n,l}$ of eigenvalues increases polynomially with respect to $\lambda_{n,l}$.
\begin{Lem}\label{lem-ev}
There exists a constant $C>0$ which satisfies
\begin{align}
\left|\lambda_{n,l\pm1}-\lambda_{n,l}\right|\geq C\lambda_{n,l}^{(c-1)/2c}
\end{align}
for any n and large $l$. 
\end{Lem}
The proof of the above lemma is found in \cite[Lemma 3.3]{Y1}, for which we use
\begin{Prop}\label{prop-weyl}
Let $d\geq3$ and $n$ be fixed.
Then as $l\to\infty$
\begin{align}\label{wkb}
\frac{1}{\pi}\int^{X_{n,l}}_{0} \left\{\lambda_{n,l}-\widetilde{V}(r)\right\}^{1/2}dr
=l+\frac{n}{2}+\frac{d}{4}+o\left(1\right),
\end{align}
where $\widetilde{V}(X_{n,l})=\lambda_{n,l}$. 
\end{Prop}
We show an outline of the proof of Proposition~\ref{prop-weyl} in the rest of this section for the readers' convenience.
Note that it has been proved for the case $d=3$, see Titchmarsh \cite[p.151]{T}.
\par
Proposition~\ref{prop-weyl} is proved by the same method as in \cite[p.158 $\S \ 7.13$]{T}.
We first observe the asymptotic behavior of solutions of \eqref{eigeneq} which are $L^2$ near $\infty$ and $0$.
%
\begin{Lem}\label{inf}
For sufficiently large $\lambda$, there is $\psi_\lambda\in H^2((1/2,\infty))$ such that
\[
\{f\in L^2((1/2,\infty)) \mid H_n f=\lambda f \text{ on } (1/2,\infty)\} = \{\alpha\psi_\lambda \mid \alpha\in\mathbb{C}\}
\]
and
\begin{align}
\psi_\lambda(1)=&2e^{-\frac{2}{3}\pi i}\left\{\lambda-U_{d,n}(1)\right\}^{-\frac{1}{4}}\left\{\cos\left(Z-\frac{1}{4}\pi\right)+O(Z^{-1})\right\}, \label{ef-from-infty}\\
\psi'_\lambda(1)=&2e^{-\frac{2}{3}\pi i}\left\{\lambda-U_{d,n}(1)\right\}^{\frac{1}{4}}\left\{\sin\left(Z-\frac{1}{4}\pi\right)+O(Z^{-1})\right\}, \label{ef'-from-infty}
\end{align}
where
\[
Z=\int_1^{T}\sqrt{\lambda-U_{d,n}(s)}ds,
\]
and $T=T_\lambda>1$ is the value satisfying $U_{d,n}(T)=\lambda$.
\end{Lem}
\begin{proof}
Let $\lambda$ be sufficiently large, and let $T=T_\lambda>1$ be as above.
We set
\begin{equation}
\zeta(r)=\int_{T}^{r}\sqrt{\lambda-U_{d,n}(s)}ds,
\label{defOfZeta}
\end{equation}
where the branch is chosen such that $\text{arg}\ \zeta(r)=\pi/2$ for $r>T$, and $\text{arg}\ \zeta(r)=-\pi$ for $r<T$.
\par
We employ Langer's method in Chapter 22.27 of \cite{T2}, where the term $\frac{5}{36\zeta^{2}}$ plays the essential role in computations of iterations near $\zeta=0$, i.e., $r=T$.
Changing variables $r\to\zeta$ and $\eta=a(r)^{-1}\psi(r)$, where $a(r)=\left(\lambda-U_{d,n}(r)\right)^{-1/4}$, implies, at least formally,
\begin{align}\label{eq infty}
\frac{d^{2}\eta}{d\zeta^{2}}+\left(1+\frac{5}{36\zeta^{2}}\right)\eta=g(\zeta)\eta ,
\end{align}
where $g(\zeta)=\frac{5}{36\zeta^{2}}-\mathcal{U}(r(\zeta),\lambda)$ and
\begin{align*}
\mathcal{U}(r,\lambda)=\frac{U^{\prime\prime}_{d,n}(r)}{4(\lambda-U_{d,n}(r))^{2}}+\frac{5(U^{\prime}_{d,n}(r))^{2}}{16(\lambda-U_{d,n}(r))^{3}} .
\end{align*}
It is well-known that the functions
\begin{align*}
\eta_1(\zeta)=\left(\frac{1}{2}\pi\zeta\right)^{\frac{1}{2}}J_{\frac{1}{3}}(\zeta), \quad 
\eta_2(\zeta)=\left(\frac{1}{2}\pi\zeta\right)^{\frac{1}{2}}H_{\frac{1}{3}}^{(1)}(\zeta),
\end{align*}
where $J_{\nu}$ is the Bessel function of the first kind and $H_{\nu}^{(j)}$ is the Hankel function, are the solutions of
\begin{align*}
\frac{d^{2}\eta}{d\zeta^{2}}+\left(1+\frac{5}{36\zeta^{2}}\right)\eta=0.
\end{align*}
We look for a solution of \eqref{eq infty} in $(\zeta(\frac12),0]\cup i[0,\infty)$ with $\eta(ix) \sim \eta_2(ix)$ as $x\to\infty$, since $\psi_\lambda(r)=a(r)\eta_2(\zeta(r))$ is $L^2$ near $\infty$.
Thanks to the standard ODE calculus, we only have to observe the solution of the integral equation
\begin{align}
\eta(\zeta) =& \eta_2(\zeta) + C_0^{-1} \int_{\zeta}^{i\infty} (\eta_1(\xi)\eta_2(\zeta) - \eta_2(\xi)\eta_1(\zeta)) g(\xi) \eta(\xi) d\xi  \nonumber\\
=& \eta_2(\zeta) + C_0^{-1} \int_{r(\zeta)}^{\infty} (\eta_1(\zeta(s))\eta_2(\zeta) - \eta_2(\zeta(s))\eta_1(\zeta)) g(\zeta(s)) \eta(\zeta(s)) \sqrt{\lambda-U_{d,n}(s)} ds,
\label{integral-eq-3-9}
\end{align}
where $C_0=W(\eta_1,\eta_2)=\eta_1(\zeta) \eta'_2(\zeta) - \eta'_1(\zeta) \eta_2(\zeta)$ is the Wronskian which is independent of $\zeta$, and the contour of the above first integral is taken over the path connecting $\zeta$ and $i\infty$  in $(-\infty,0]\cup i[0,\infty)$.

We write the second term of \eqref{integral-eq-3-9} by $F_\lambda \eta(\zeta)$, and we define a function space $X$ and its norm by
\begin{align*}
\|\eta\|_{X}=\sup_{\zeta\in(\zeta(\frac12),0]\cup i[0,\infty)} e^{\Im\zeta} |\eta(\zeta)|,\quad
X=\left\{\eta\in C\left((\zeta(1/2),0]\cup i[0,\infty);\mathbb C\right):\|\eta\|_{X}<\infty\right\}.
\end{align*}
We note that $\eta_1(\zeta)=O(e^{\text{Im}\zeta})$, $\eta_2(\zeta)=O(e^{-\text{Im}\zeta})$ (see \cite[\S 7.21]{W} for the asymptotic behavior of $J_{\frac{1}{3}}(\zeta)$ and $H_{\frac{1}{3}}^{(1)}(\zeta)$ as $|\zeta|\to\infty$), and that
\begin{align}
\int_{1/2}^\infty |g(\zeta(r))| |\lambda-U_{d,n}(r)|^{\frac{1}{2}} dr =O(\lambda^{-1/2}T^{-1}) \quad (\lambda\to\infty)\label{eq-appendix}
\end{align}
(see Appendix for the proof).
Then we have $\|\eta_2\|_{X}<\infty$ and $\|F_\lambda\|_{X\to X}=O(\lambda^{-1/2}T^{-1})$, which implies
\begin{align}\label{eta esti}
\eta(\zeta) = \eta_2(\zeta) + O(\lambda^{-\frac{1}{2}}T^{-1}e^{-\text{Im}\zeta})
\end{align}
for sufficiently large $\lambda$ and $\zeta\in(\zeta(\frac12),0]\cup i[0,\infty)$, and 
thus
\begin{align*}
\psi_\lambda(r)=\left\{\lambda-U_{d,n}(r)\right\}^{-\frac{1}{4}}\left\{ \left(\frac{1}{2}\pi\zeta\right)^{\frac{1}{2}}H^{(1)}_{\frac{1}{3}}(\zeta) +O(\lambda^{-\frac{1}{2}}T^{-1}e^{-\text{Im}\zeta})\right\}.
\end{align*}


If $r<T$, we have $H_{\frac{1}{3}}^{(1)}(\zeta)=\frac{2}{\sqrt{3}}e^{-\frac{1}{6}\pi i}\{J_{\frac{1}{3}}(z)+J_{-\frac{1}{3}}(z)\}$, where $\zeta=e^{-i\pi}z$ (see \cite[\S 7.8]{T}). 
By using 
\begin{align}\label{Bessel infty}
J_{\nu}(z)=\left(\frac{2}{\pi z}\right)^{\frac{1}{2}}\left\{\cos\left(z-\frac{1}{2}\nu-\frac{1}{4}\pi\right)+O(z^{-1})\right\}\quad \text{as } z\to\infty, 
\end{align}
(see \cite[\S 7.21]{W}), we have 
\begin{align*}
J_{\frac{1}{3}}(z)+J_{-\frac{1}{3}}(z)=\left(\frac{6}{\pi z}\right)^{\frac{1}{2}}\left\{\cos\left(z-\frac{1}{4}\pi\right)+O(z^{-1})\right\}. 
\end{align*}
Since $Z=\zeta(1)=\int_1^{T}\sqrt{\lambda-U_{d,n}(s)}ds\leq\lambda^{1/2} T$ for sufficiently large $\lambda$, we obtain \eqref{ef-from-infty}.

For the proof of \eqref{ef'-from-infty}, differentiating \eqref{integral-eq-3-9} in $\zeta$ on $(\zeta(\frac12),0)\cup i(0,\infty)$ implies
\begin{align*}
\eta'(\zeta)
= \eta'_2(\zeta) + C_0^{-1} \int_{r(\zeta)}^{\infty} (\eta_1(\zeta(s))\eta'_2(\zeta) - \eta_2(\zeta(s))\eta'_1(\zeta)) g(\zeta(s)) \eta(\zeta(s)) \sqrt{\lambda-U_{d,n}(s)} ds.
\end{align*}
Thus, by \eqref{eta esti}, we have
\begin{align}\label{eta' esti}
\eta'(\zeta)=(1-O(\lambda^{-\frac{1}{2}}T^{-1}))\eta'_2(\zeta)+O(\lambda^{-\frac{1}{2}}T^{-1})\eta'_1(\zeta),
\end{align}
and therefore the identity $J_\nu'(z)=-J_{\nu+1}(z) + \frac{\nu}{z}J_\nu(z)$, \eqref{Bessel infty} and $\frac{d\zeta}{dr}=\sqrt{\lambda-U_{d,n}(r)}$ imply \eqref{ef'-from-infty}.

\end{proof}


%
\begin{Lem}\label{zero}
For sufficiently large $\lambda$, there is $\phi_\lambda\in H^2((0,2))$ such that
\begin{align*}
&\{f\in L^2((0,2)) \mid H_n f=\lambda f \text{ on } (0,2)\} = \{\alpha\phi_\lambda \mid \alpha\in\mathbb{C}\} , \quad (d,n)\neq(3,0), \\
&\{f\in H^1((0,2)) \mid H_n f=\lambda f \text{ on } (0,2), \  f(0)=0\} = \{\alpha\phi_\lambda \mid \alpha\in\mathbb{C}\}, \quad (d,n)=(3,0),
\end{align*}
and
\begin{align}
\phi_\lambda(1)=&\sqrt{\frac{2}{\pi}}\lambda^{-\frac{1}{4}}\left\{\cos\left(\sqrt{\lambda}-\frac{2n+d-1}{4}\pi\right)
+O\left(\lambda^{-\frac{1}{2}}\right)\right\}, \label{ef-from-0} \\
\phi'_\lambda(1)=&\sqrt{\frac{2}{\pi}}\lambda^{\frac{1}{4}}\left\{-\sin\left(\sqrt{\lambda}-\frac{2n+d-1}{4}\pi\right)
+O\left(\lambda^{-\frac{1}{2}}\right)\right\}. \label{ef'-from-0}
\end{align}
\end{Lem}
\begin{proof}
We suppose $(d,n)\neq(3,0)$ for simplicity.
It is easy to see that $\phi_+(r)=r^{\frac{1}{2}}J_{(n+\frac{d-2}{2})}(r\sqrt{\lambda})$ and $\phi_{-}(r)=r^{\frac{1}{2}}Y_{(n+\frac{d-2}{2})}(r\sqrt{\lambda})$, where $Y_{\nu}$ is the Bessel function of the second kind, are solutions of
\begin{align*}
\frac{d^{2}}{dr^{2}}f(r)+\left\{\lambda - \frac{(d-1)(d-3)}{4r^{2}} - \frac{n(n+d-2)}{r^{2}}\right\}f(r)= 0,
\end{align*}
and
\[
\phi_+(r)=r^{n+\frac{d-1}{2}}(C'_{\lambda}+o(1)),
\quad
\phi_-(r)=
r^{-n-\frac{d-3}{2}}(-\frac{1}{\pi C'_{\lambda}}+o(1)) \quad \text{as }r\to0,
\]
where $C'_{\lambda} = \frac{1}{\Gamma(n+\frac{d}{2})} \left(\frac{\sqrt{\lambda}}{2}\right)^{n+\frac{d-2}{2}}$ (see \cite[\S 3]{W}).
Thus it suffices to construct the solution $\phi$ of \eqref{eigeneq} with $\phi\sim\phi_+$, regarding the term $\tilde{V}(r)\phi(r)$ as a perturbation employing the norm $\|f\|=\sup_{r\in(0,2)}\min(1,(\sqrt{\lambda} r)^{n+\frac{d-1}{2}})^{-1}|f(r)|$.
Finally we have \eqref{ef-from-0} and \eqref{ef'-from-0} by the asymptotic behavior \eqref{Bessel infty}.
\end{proof}
\begin{proof}[Proof of Proposition \ref{prop-weyl}]
Suppose that $\lambda_{n,l}$ is an eigenvalue of $H_n$ and that $f_{n,l}$ is the eigenfunction for $\lambda_{n,l}$. 
Since Lemmas~\ref{inf} and \ref{zero} imply $f_{n,l}=A\psi_{\lambda_{n,l}}=B\phi_{\lambda_{n,l}}$ on $(1/2, 2)$ with some constants $A$ and $B$, we have
\begin{align*}
\phi_{\lambda_{n,l}}(1)\psi^{\prime}_{\lambda_{n,l}}(1)-\phi_{\lambda_{n,l}}^{\prime}(1)\psi_{\lambda_{n,l}}(1) = 0 .
\end{align*}
By \eqref{ef-from-infty}, \eqref{ef'-from-infty}, \eqref{ef-from-0} and \eqref{ef'-from-0}, we have
\begin{align*}
\sin\left(Z+\sqrt{\lambda_{n,l}}-\frac{2n+d}{4}\pi\right)
=O\left(\lambda_{n,l}^{-\frac{1}{2}}\right)+O(Z^{-1}).
\end{align*}
If we write
\begin{align*}
Z+\sqrt{\lambda_{n,l}}=\left(\frac{2n+d}{4}+m_{n,l}\right)\pi+\delta
\end{align*}
with some $m_{n,l}\in\mathbb{Z}$ and $\delta\in[-\frac{\pi}{2},\frac{\pi}{2})$, then
\begin{align*}
\sin\delta
=O\left(\lambda_{n,l}^{-\frac{1}{2}}\right)+O(Z^{-1})
\end{align*}
and so 
\begin{align*}
\delta=
O\left(\lambda_{n,l}^{-\frac{1}{2}}\right)+O(Z^{-1}).
\end{align*}
Then we have 
\begin{align}\label{Z1}
Z+\sqrt{\lambda_{n,l}}=\left(\frac{2n+d}{4}+m_{n,l}\right)\pi
+O\left(\lambda_{n,l}^{-\frac{1}{2}}\right)+O(Z^{-1}).
\end{align}
On the other hand, it is proved in \cite[p161-163]{T} that, letting $\widetilde{V}(X_{n,l})=\lambda_{n,l}$,
\begin{align}\label{Z2}
Z+\sqrt{\lambda_{n,l}}
=\int_{0}^{X_{n,l}}\{\lambda_{n,l}-\widetilde{V}(r)\}^{\frac12}dr + O(\lambda_{n,l}^{-\frac12}). 
\end{align}
Thus, using \eqref{Z1} and \eqref{Z2}, we obtain 
\[
\int_{0}^{X_{n,l}}\{\lambda_{n,l}-\widetilde{V}(r)\}^{\frac12}dr
=\left(\frac{2n+d}{4}+m_{n,l}\right)\pi
+O\left(\lambda_{n,l}^{-\frac{1}{2}}\right).
\]

It suffices to show that $m_{n,l}=l$ for sufficiently large $l$, which follows from counting the zeros of $f_{n,l}$ in $(0,\infty )$ by two different methods. 

We easily see (e.g.\ \cite[\S 5.4]{T}) that $f_{n,l}$ has $l$ zeros in $(0,\infty )$, $f_{n,l}$ is the $l+1$-th eigenfunction of $H_n$.
On the other hand, we can find that $m_{n,l}$ is the number of the zeros of $f_{n,l}$
, which concludes the proof.
In the rest of the proof, we show that
\begin{itemize}
\item $f_{n,l}$ has $m_{n,l}-p$ zeros in $[r_{0},T)$,
\item $f_{n,l}(r)$ has $p$ zeros in $(0,r_0)$,
\item $f_{n,l}$ has no zeros in $[T,\infty)$,
\end{itemize}
where $p=\left[\pi^{-1}\sqrt{\lambda_{n,l}}-\frac{2n+d-1}{4}\right]$ and $r_0=(p+\frac{2n+d-1}{4})\pi/\sqrt{\lambda_{n,l}}$.
We note that $1-\frac{1}{\pi\sqrt{\lambda_{n,l}}}<r_0\leq1$.

The third assertion is proved by contradiction. If there exists a zero of $f_{n,l}$ in $[T,\infty)$, $f_{n,l}$ and $f'_{n,l}$ are positive (or negative) near the zero. Since $f''_{n,l}(r)=(U_{d,n}(r)-\lambda)f_{n,l}(r)$, $f_{n,l}$ is convex (or concave). Hence $f_{n,l}$ tends to $\infty$ (or $-\infty$). This contradicts that $f_{n,l}\in L^{2}((0,\infty))$. 

For the first assertion, we remark that $f_{n,l}(r)=0$ if and only if $\eta(\zeta(r))=a(r)^{-1}f_{n,l}(r)=0$, where $a(r)$ and $\zeta(r)$ are as in the proof of Lemma~\ref{inf} and $\zeta([r_0,T))=[\zeta(r_0),0)$.

We observe that the number of the zeros of $J_{\frac{1}{3}}(r)+J_{-\frac{1}{3}}(r)$ is $m_{n,l}-p$ in $(0,\left(m_{n,l}-p+\frac{1}{4}\right)\pi]$ for large $l$ (\cite[Lemma 7.9b]{T}), and $-\zeta(r_{0})=\left(m_{n,l}-p+\frac{1}{4}\right)\pi+o(1)$.
Indeed, since
\begin{align*}
-\zeta (r_{0})-Z=\int_{r_0}^{1}\{\lambda_{n,l}-U_{d,n}(r)\}^{\frac12}dr
=(1-r_{0})\sqrt{\lambda_{n,l}}+O(\lambda_{n,l}^{-\frac{1}{2}}),
\end{align*}
the definition of $m_{n,l}$ shows that 
\begin{align*}
-\zeta(r_{0})=\left(\frac{2n+d}{4}+m_{n,l}\right)\pi-r_{0}\sqrt{\lambda_{n,l}}+o(1)=\left(m_{n,l}-p+\frac{1}{4}\right)\pi+o(1).
\end{align*}
Thus $\eta_2(\zeta)=\left(\frac{1}{2}\pi\zeta\right)^{\frac{1}{2}}H_{\frac{1}{3}}^{(1)}(\zeta)$ has exactly $m_{n,l}-p$ zeros in $[\zeta(r_0),0)$.
Moreover, $\eta_2$ is concave (resp.\ convex) if $\eta_2(\zeta)>0$ (resp.\ $\eta_2(\zeta)<0$) since $\eta_2$ solves $\frac{d^{2}\eta}{d\zeta^{2}}+\left(1+\frac{5}{36\zeta^{2}}\right)\eta=0$.
Hence the asymptotic behavior of $\eta_2(\zeta)=C\left\{\cos\left(-\zeta-\frac{1}{4}\pi\right)+O(|\zeta|^{-1})\right\}$  implies
\[
\ep=\inf\{\text{maximal values of }|\eta_2(\zeta)| \text{ in }(-\infty,0)\}>0.
\]

Let $S=\{\zeta\in[\zeta(r_0),0) \mid |\eta(\zeta)| < \ep/3\}$.
Then $\eta(\zeta)$ is monotone on $\zeta\in S$.
In fact, $S\subset \eta_2^{-1}((-2\ep/3,2\ep/3))$ for large $l$ by \eqref{eta esti}, and thus \eqref{eta' esti} implies $\eta'(\zeta)\neq0$ for $\zeta\in S$.
Therefore each segment of $S$, except the one of the form $(-a,0)$, has exactly one zero, and the number of segments of $S\setminus(-a,0)$ is $m_{n,l}-p$, same as that of zeros of $\eta_2$ thanks to \eqref{eta esti}.

The second assertion is proved by the same argument as above, noting that $\phi_+(r)=r^{\frac{1}{2}}J_{(n+\frac{d-2}{2})}(r\sqrt{\lambda})$ has $p$ zeros in $(0,r_0)$ (the number of the zeros of $J_{(n+\frac{d-2}{2})}(r)$ is $p$ in $(0,(p+\frac{2n+d-1}{4})\pi)$ by \cite[Lemma~7.9a]{T}).
\end{proof}
%
%
%
\section{Proof of Theorem~\ref{main} and Proposition~\ref{cor-main}}\label{sec-proof}
%
We first note the formula for $\kappa \in C_{0}^{\infty}(\R)$ and $\Phi, \Psi \in C_{0}^{\infty}(\R^d)$
\begin{align}\label{eq-formula-FMS}
\int_{\R} \int_{\R^d} \int_{\R^d} E(t,x,y) \kappa(t) \Phi(x) \Psi(y) dy dx dt
= (\hat\kappa(H) \Psi , \overline{\Phi})_{L^2(\R^d)}
= (\hat\kappa(H) \Phi , \overline{\Psi})_{L^2(\R^d)}
\end{align}
with the help of identity $U(t)=e^{-itH}$, the spectral decomposition theorem and Fubini's theorem, where $\hat{\kappa}(\lambda)=\int_{\R} e^{-it\lambda}\kappa(t) dt$.
It follows from Lemma~\ref{lem-decomp} that
\begin{align}\label{eq-formula-FMS2}
\int_{\R} \int_{\R^d} \int_{\R^d} E^{Y_{nm}}(t,x,y) \kappa(t) \Phi(x) \Psi(y) dy dx dt
= (\hat\kappa(H_n)J_{nm}\Phi , J_{nm}\overline{\Psi})_{L^2((0,\infty))} .
\end{align}
\begin{proof}[Proof of Proposition~\ref{cor-main}]
It suffices to show that $E^{Y_n}(t,x,y)$ is not in $C^1$ on $(0,\infty)\times\{x\in\mathbb R^{d}\mid x\neq0\}\times\{y\in\mathbb R^{d}\mid y\neq 0\}$, which we show by contradiction.
Suppose that $E^{Y_n}$ is in $C^1$ near $(t_0,x_0,y_0)$ with $t_0>0,\ x_0 \neq0,\ \text{and}\ y_0 \neq 0$.
The definition of $E^{Y_n}$ and Remark~\ref{rem-thm} imply that 
\begin{align}\label{C^1}
E^{Y_n}\in C^1(I \times B_{r_1,R_1} \times B_{r_2,R_2})
\end{align}
with some interval $I \Subset (0,\infty)$ and $0 < r_j<R_j$, $j=1,2$, where $B_{r,R} = \{ x\in\R^d \mid r<|x|<R\}$.
We may assume $Y_n=Y_{n1}$ without loss of generality.
\par
We fix $\kappa \in C_{0}^{\infty}(I)$, $\Phi \in C_{0}^{\infty}(B_{r_1,R_1})$ and $\Psi \in C_{0}^{\infty}(B_{r_2,R_2})$ so that $\hat\kappa(0)\neq0$ and that $\phi(r)=J_{n1}\Phi(r) \not\equiv0$ and $\psi(r)=J_{n1}\Psi(r) \not\equiv0$ are non-negative.
We set
\begin{align*}
G(\tau,j,k)=&\int_{\mathbb R} \int_{\mathbb R^{d}} \int_{\mathbb R^{d}}E^{Y_{n}}(t,x,y)\kappa(t)\Phi(x)\Psi(y)e^{it\tau+i|x|j-i|y|k}dydxdt \\
=& (\hat\kappa(H_n-\tau)\phi_j , \overline{\psi_k})_{L^2((0,\infty))}
\end{align*}
for $(\tau,j,k) \in \mathbb R \times \mathbb R \times \mathbb R$, where we used \eqref{eq-formula-FMS2} and the notation
\begin{align*}
\phi_j(r)= e^{irj}\phi(r), \quad \psi_k(r)= e^{-irk}\psi(r).
\end{align*}
It suffices to show that there exist a constant $C>0$ and a sequence $(\tau_{l},j_{l},k_{l})$ such that $\left|\tau_{l}\right|+\left|j_{l}\right|+\left|k_{l}\right| \rightarrow \infty$, as $l\rightarrow\infty$ such that
\begin{align}\label{eq-contra}
|G(\tau_{l},j_{l},k_{l})|\geq C\left(1+|\tau_{l}|+|j_{l}|+|k_{l}|\right)^{-1/2c},
\end{align}
which contradicts \eqref{C^1}.
In the following, we show \eqref{eq-contra} by taking $\tau=\lambda_{n,l}$ and $j_l=k_l=\sqrt{\lambda_{n,l}}$.
\par
The spectral decomposition theorem implies
\begin{align*}
\hat\kappa(H_n-\tau)u = \sum_{l=0}^\infty \hat{\kappa}(\lambda_{n,l}-\tau) (u, f_{n,l})_{L^2((0,\infty))} f_{n,l}, \quad u \in L^2((0,\infty))
\end{align*}
where $f_{n,l}$, $l=0,1,\dots$, is the normalized eigenfunction of $H_n$ with eigenvalue $\lambda_{n,l}$.
Then it follows from Lemma~\ref{lem-ev} and $\widehat{\kappa}\in\mathcal{S}(\R)$ that for any $N\in\mathbb N$
\begin{align*}
&\| \hat\kappa(H_n-\tau_l) -\hat{\kappa}(0) (\cdot, f_{n,l})_{L^2((0,\infty))} f_{n,l} \|_{L^2((0,\infty)) \to L^2((0,\infty))} \\
=&\sup_{l'\in\mathbb{N}_+ \setminus \{l\}} |\hat{\kappa}(\lambda_{n,l'}-\lambda_{n,l})| \\
\leq& C_N \lambda_{n,l}^{-N},
\end{align*}
which implies
\begin{align*}
G(\tau_l,j_l,k_l)
= \hat{\kappa}(0) (\phi_{j_l}, f_{n,l})_{L^2((0,\infty))} (\psi_{k_l}, \overline{f_{n,l}})_{L^2((0,\infty))} + O(\lambda_{n,l}^{-\infty}).
\end{align*}
We note that $\operatorname{supp} \phi \Subset \Omega_{\lambda}$ for sufficiently large $\lambda$, where $\Omega_{\lambda}$ is as in Lemma~\ref{lem-cls-allowed}.
Thus we can compute
\begin{align*}
(\phi_{j_l}, f_{n,l})_{L^2((0,\infty))}
=& \int_0^\infty e^{irj_l} \phi(r) \cdot \lambda_{n,l}^{-1/4}\left(\Re\{C_{\lambda_{n,l}}e^{i\sqrt{\lambda_{n,l}}r}\}+O(C_{\lambda_{n,l}}\lambda_{n,l}^{-1/2})\right) dr \\
=& \frac{1}{2\lambda_{n,l}^{1/4}} \left\{\int_0^\infty \left(C_{\lambda_{n,l}}e^{2i\sqrt{\lambda_{n,l}}r} + \overline{C_{\lambda_{n,l}}} \right) \phi(r) dr +O(C_{\lambda_{n,l}}\lambda_{n,l}^{-1/2}) \right\} \\
=& \frac{\overline{C_{\lambda_{n,l}}}}{2\lambda_{n,l}^{1/4}} \left( \int_0^\infty \phi(r) dr + O(\lambda_{n,l}^{-1/2}) \right),
\end{align*}
where we used integration by parts once to show $\int_0^\infty e^{2i\sqrt{\lambda_{n,l}}r} \phi(r) dr = O(\lambda_{n,l}^{-1/2})$.
The same computations imply $(\psi_{k_l}, f_{n,l})_{L^2((0,\infty))} = \frac{C_{\lambda_{n,l}}}{2\lambda_{n,l}^{1/4}} \left( \int_0^\infty \psi(r) dr + O(\lambda_{n,l}^{-1/2}) \right)$.
Thus we obtain, taking Lemma~\ref{lem-cls-allowed-const} into account,
\begin{align*}
G(\tau_l,j_l,k_l)
=& \hat{\kappa}(0) (\phi_{j_l}, f_{n,l})_{L^2((0,\infty))} (\psi_{k_l}, f_{n,l})_{L^2((0,\infty))} + O(\lambda_{n,l}^{-\infty}) \\
=& \frac{\hat{\kappa}(0) |C_{\lambda_{n,l}}|^2}{4\lambda_{n,l}^{1/2}} \left( \int_0^\infty \phi(r) dr \int_0^\infty \psi(r) dr + O(\lambda_{n,l}^{-1/2}) \right)+ O(\lambda_{n,l}^{-\infty}) \\
\geq& c_0 \lambda_{n,l}^{-1/2c}
\end{align*}
for sufficiently large $l$ with some $c_0>0$, which concludes the assertion.
\end{proof}
\begin{proof}[Proof of Theorem~\ref{main}]
Let
\begin{align*}
G_0(\tau,j,k)=&\int_{\mathbb R} \int_{\mathbb R^{d}} \int_{\mathbb R^{d}}E(t,x,y)\kappa(t)J_{00}^*\phi(x) J_{00}^*\psi(y)e^{it\tau+i|x|j-i|y|k}dydxdt ,
\end{align*}
where $\kappa \in C_0^\infty(\R)$ and $\phi,\psi \in C_0^\infty((0,\infty))$ are arbitrary non-negative functions
and $J_{00}: L^2(\mathbb{R}^d)\to L^2((0,\infty))$ is as in Section~\ref{sec-polar}.
By \eqref{eq-formula-FMS} and Lemma~\ref{lem-decomp} we have
\begin{align*}
G_0(\tau,j,k) = (\hat\kappa(H-\tau) J_{00}^*\phi_j , \overline{J_{00}^*\psi_k})_{L^2(\R^d)}
=  (\hat\kappa(H_{0}-\tau)\phi_j , \overline{\psi_k})_{L^2((0,\infty))},
\end{align*}
where $\phi_j(r)= e^{irj}\phi(r)$ and $\psi_k(r)= e^{-irk}\psi(r)$.
Thus the same argument as in the proof of Theorem~\ref{cor-main} implies $G(\tau_l,j_l,k_l)\geq c_0 \lambda_{n,l}^{-1/2c}$.
Therefore we obtain $E(t,x,y)\notin C^1(I \times B_{r_1,R_1} \times B_{r_2,R_2})$ for any interval  $I\Subset (0,\infty)$ and $0<r_j<R_j$, $j=1,2$, which completes the proof.
\end{proof}
%
%
%
\section*{Acknowledgement}
The authors would like to thank the referee for helpful comments. 
The authors would also like to thank Fumihito Abe, Ryo Muramatsu and Kouichi Taira for helpful discussions. 
KK was partially supported by JSPS KAKENHI Grant Number 22K03394.
YT was partially supported by JSPS KAKENHI Grant Number 23K12991.
The authors appreciate useful comments made by the referee.
%

\section*{Appendix: Proof of \eqref{eq-appendix}}

We write $U(r)$ instead of $U_{d,n}(r)$ until the end of the proof. 
We divide the integral as follows,
\[
\int_{1/2}^\infty |g(\zeta(r))| |\lambda-U(r)|^{\frac{1}{2}} dr = \int_{1/2}^{(1-\varepsilon)T}+\int_{(1-\varepsilon)T}^{(1+\varepsilon)T}+\int_{(1+\varepsilon)T}^{\infty}=I_{1}+I_{2}+I_{3},
\]
where $\varepsilon>0$ is a small constant which will be fixed later.\\
\underline{Estimate of $I_{2}$}

Integrating by part twice on the definition \eqref{defOfZeta} of $\zeta$, we have 
\begin{align*}
 -\zeta (r) &= -\int_{T}^{r}\sqrt{\lambda - U(s)}ds\\
&= \frac{2}{3}\int_{T}^{r}({(\lambda - U(s))}^{3/2})'\frac{1}{U'(s)}ds\\
&=\frac{2}{3}\frac{(\lambda - U(r))^{3/2}}{U'(r)}+\frac{2}{3}\int_{T}^{r}{(\lambda - U(s))}^{3/2}\frac{U''(s)}{U'(s)^{2}}ds\\
&= \frac{2}{3}\frac{(\lambda - U(r))^{3/2}}{U'(r)}-\frac{4}{15}\int_{T}^{r}({(\lambda - U(s))}^{5/2})'\frac{U''(s)}{U'(s)^{3}}ds\\
&=\frac{2}{3}\frac{(\lambda - U(r))^{3/2}}{U'(r)}\left\{1- \frac{2}{5}\frac{U''(r)(\lambda - U(r))}{U'(r)^2} + S(r) \right\},
\end{align*}
where
$$
S(r)=\frac{2U'(r)}{5(\lambda - U(r))^{3/2}}\int_{T}^{r}(\lambda -U(s))^{5/2}\frac{d}{ds}\left(\frac{U''(s)}{U'(s)^{3}}\right)ds. 
$$

We claim that for $T/2\leq r\leq 2T$,
\begin{align*}
\frac{U''(r)(\lambda - U(r))}{U'(r)^2} = O\left(\frac{\lambda-U(r)}{U(r)}\right) , \quad
S(r) = O\left(\left(\frac{\lambda-U(r)}{U(r)}\right)^2\right) .
\end{align*}
The first assertion is easy to show by \eqref{Ass2} and \eqref{Ass3}.
For the second one, \eqref{Ass3} implies
\begin{align*}
\frac{d}{ds}\left(\frac{U''(s)}{U'(s)^{3}}\right)
=\frac{U'''(s)}{U'(s)^{3}}-\frac{3U''(s)^{2}}{U'(s)^{4}}
=\frac{1}{U'(s)^{2}}O(s^{-2}).
\end{align*}
Since $U'(r)$ is monotonically increasing, we have, for any $T\leq r\leq 2T$, 
\begin{align*}
|S(r)|&=\frac{2U'(r)}{5(U(r)-\lambda)^{3/2}}\int_{T}^{r}(U(s)-\lambda)^{5/2} \left|\frac{d}{ds}\left(\frac{U''(s)}{U'(s)^{3}}\right)\right| ds\\
&\leq C\frac{U'(r)}{(U(r)- \lambda)^{3/2}}\int_{T}^{r} \frac{(U(s)-\lambda)^{5/2}}{s^2U'(s)^2} ds\\
&= \frac{2C}{7}\cdot\frac{U'(r)}{(U(r)- \lambda)^{3/2}}\int_{T}^{r} \frac{\left((U(s)-\lambda)^{7/2}\right)'}{s^2U'(s)^3} ds\\
&\leq C'\frac{U'(r)}{(U(r)- \lambda)^{3/2}}\cdot\frac{1}{T^{2}U'(T)^{3}}\int_{T}^{r}\left((U(s)-\lambda)^{7/2}\right)^{'}ds\\
&\leq C''\frac{U'(r)(\lambda - U(r))^{2}}{U(T)^{2}U'(T)}\\
&=C'' \left(\frac{\lambda - U(r)}{U(r)}\right)^2  \left(\frac{U(r)}{U(T)}\right)^2 \frac{U'(r)}{U'(T)}.
\end{align*}
By the same argument, we have, for any $T/2\leq r\leq T$,
\begin{align*}
|S(r)|\leq C \left(\frac{\lambda - U(r)}{U(r)}\right)^2 .
\end{align*}
Since \eqref{Ass3} implies $\frac{U(2r)}{U(r)}=O(1)$ and $\frac{U'(2r)}{U'(r)}=O(1)$, we have the desired estimate.

We can fix $\ep>0$ so that
\begin{align*}
\left|\frac{2}{5}\frac{U''(r)(\lambda - U(r))}{U'(r)^2} - S(r)\right|\leq\frac{1}{2}
\end{align*}
for $(1-\varepsilon)T\le r \le (1+\varepsilon)T$.
In fact, $\log\frac{U(ar)}{U(r)}= \int_r^{ar} \frac{U'(s)}{U(s)}ds$ and \eqref{Ass2}, \eqref{Ass3} imply for $a\geq1$
\begin{align*}
0\leq\log\frac{U(aT)}{U(T)} \leq C \int_r^{ar} \frac{ds}{s} = C\log a.
\end{align*}
Thus we can make $\left|\frac{\lambda-U(r)}{U(r)}\right| = \left|1-\frac{U(T)}{U(r)}\right|$ sufficiently small as $\ep\to0$. 
By using Taylor's expansion of $(1-x)^{-2}$, we have
\begin{align*}
\frac{1}{\zeta (r)^{2}} &= \frac{9U'(r)^{2}}{4(\lambda - U(r))^{3}}\left\{1+ \frac{4}{5}\frac{U''(r)(\lambda - U(r))}{U'(r)^2} + O\left(\frac{(\lambda - U(r))^{2}}{U(r)^2}\right)\right\}\\
&=\frac{9U'(r)^{2}}{4(\lambda - U(r))^{3}}+\frac{9U''(r)}{5(\lambda - U(r))^{2}}+O\left(T^{-2}|\lambda -U(r)|^{-1}\right).
\end{align*}
Hence we obtain
\begin{align*}
\frac{5}{36\zeta (r)^{2}} = \mathcal{U}(r, \lambda ) + 
O\left(T^{-2}|\lambda -U(r)|^{-1}\right),
\end{align*}
which implies that $I_{2}=O(\lambda^{-1/2}T^{-1})$. \\
\underline{Estimate of $I_{1}$}

Let $R$ be the constant as in Assumption~\ref{Ass}.
By \eqref{Ass2}, there exists a constant $C>0$ such that $\lambda=U(T)\geq CT^{2c}$, i.e., $\lambda^{-\frac{1}{2c}}\leq CT^{-1}$. 
Then we have 
\begin{align*}
\int_{1/2}^{R}|\mathcal U(t,\lambda)||\lambda-U(t)|^{1/2}dt
&\leq\int_{1/2}^{R}\frac{|U''(t)|}{4(\lambda-U(t))^{3/2}}+\frac{5U'(t)^{2}}{16(\lambda-U(t))^{5/2}}dt\\
&\leq C\lambda^{-3/2}\\
&\leq C\lambda^{-1/2}T^{-1}.
\end{align*}
Integrating by part yields that 
\begin{align*}
\int_{R}^{(1-\varepsilon)T}\frac{U''(t)}{(\lambda-U(t))^{3/2}}dt
&=\left[\frac{U'(t)}{(\lambda-U(t))^{3/2}}\right]^{(1-\varepsilon)T}_{R}-\frac{3}{2}\int_{R}^{(1-\varepsilon)T}\frac{U'(t)^{2}}{(\lambda-U(t))^{5/2}}dt\\
&=\frac{U'((1-\varepsilon)T)}{(\lambda-U((1-\varepsilon)T))^{3/2}}-\frac{U'(R)}{(\lambda-U(R))^{3/2}}-\frac{3}{2}\int_{R}^{(1-\varepsilon)T}\frac{U'(t)^{2}}{(\lambda-U(t))^{5/2}}dt,
\end{align*}
which implies that 
\begin{align*}
\int_{R}^{(1-\varepsilon)T}\mathcal U(t,\lambda)&\{\lambda-U(t)\}^{1/2}dt\\
&=\frac{U'((1-\varepsilon)T)}{4(\lambda-U((1-\varepsilon)T))^{3/2}}-\frac{U'(R)}{4(\lambda-U(R))^{3/2}}-\frac{1}{16}\int_{R}^{(1-\varepsilon)T}\frac{U'(t)^{2}}{(\lambda-U(t))^{5/2}}dt\\
&\leq\frac{U'((1-\varepsilon)T)}{2(\lambda-U((1-\varepsilon)T))^{3/2}}
\leq C_{\varepsilon}\lambda^{-1/2}T^{-1}.
\end{align*}
By the definition \eqref{defOfZeta} of $\zeta$, we have 
\begin{align*}
\int_{1/2}^{(1-\varepsilon)T}\zeta(t)^{-2}\{\lambda-U(t)\}^{1/2}dt=\zeta((1-\varepsilon)T)^{-1}-\zeta(1/2)^{-1}
\leq-\zeta(1/2)^{-1}
\leq C\lambda^{-1/2}T^{-1}. 
\end{align*}
Hence we have $I_{1}=O(\lambda^{-1/2}T^{-1})$.
\\
\underline{Estimate of $I_{3}$}\\
Since \eqref{Ass2} yields that $r^{-1}U(r)$ is monotonically increasing, we have
\begin{align*}
|\zeta((1+\varepsilon)T)|&=\int_{T}^{(1+\varepsilon)T}\sqrt{U(s)-\lambda}ds\\
&\geq\int_{T}^{(1+\varepsilon)T}\left\{(s-T)\left(\min_{r\in[T,(1+\varepsilon)T]}U'(r)\right)\right\}^{1/2}ds\\
&\geq\left\{\frac{U(T)}{T}\right\}^{1/2}\int_{T}^{(1+\varepsilon)T}(s-T)^{1/2}ds= \frac23 \varepsilon^{3/2}\lambda^{1/2}T, 
\end{align*}
which shows that
$$
\int_{(1+\varepsilon)T}^{\infty} \frac{(\lambda-U(s))^{1/2}}{\zeta(s)^2} ds=\frac{1}{\zeta((1+\varepsilon)T)}=O(\lambda^{-1/2}T^{-1}).
$$
Noting that \eqref{Ass2} and $\log\frac{U(ar)}{U(r)}= \int_r^{ar} \frac{U'(s)}{U(s)}ds$ imply $U((1+\ep)T)\geq c_\ep U(T)=c_\ep \lambda$ with some $c_\ep>1$, we have
\begin{align*}
\int_{(1+\varepsilon)T}^{\infty}\frac{(U'(s))^{2}}{(U(s)-\lambda)^{5/2}}ds
&=\int_{(1+\varepsilon)T}^{\infty}\left(\frac{U'(s)}{U(s)}\right)^{2}\times\left(\frac{U(s)}{U(s)-\lambda}\right)^{2}\times\frac{1}{(U(s)-\lambda)^{1/2}}ds\\
&\leq C_{\varepsilon}\lambda^{-1/2}\int_{(1+\varepsilon)T}^{\infty}s^{-2}ds\leq C_{\varepsilon}\lambda^{-1/2}T^{-1}.
\end{align*}
We also have the following estimate
$$
\int_{(1+\varepsilon)T}^{\infty}\frac{U''(s)}{(U(s)-\lambda)^{3/2}}ds
=\int_{(1+\varepsilon)T}^{\infty}\frac{U''(s)}{U'(s)}\frac{U'(s)}{(U(s)-\lambda)^{3/2}}ds
\leq C_{\varepsilon}\lambda^{-1/2}T^{-1}.
$$
Hence we obtain $I_{3}=O(\lambda^{-1/2}T^{-1})$.
\qed


\noindent
\end{document}